\documentclass[11pt]{article}
\usepackage{amsfonts, amsmath,amssymb}
\usepackage{verbatim}
\usepackage[dvips]{graphicx}

\newcommand{\detail}[1]{\par\noi{\bf [Proof detail\ }{#1}
\hfill{\bf ]}\par\noi\hspace{-4pt}}
\renewcommand{\detail}[1]{}

\newcommand{\file}{$\ti{\ }$/wisk/anja/IF$\underline{\ }$cor.tex\quad}
\renewcommand{\file}{}

\newcommand{\dis}{\displaystyle}
\newcommand{\txt}{\textstyle}

\newcommand{\noi}{\noindent}

\newcommand{\halmos}{\rule{1ex}{1.4ex}}
\def \qed {\nopagebreak{\hspace*{\fill}$\halmos$\medskip}}
\newcommand{\med}{\medskip}

\newtheorem{theorem}{Theorem}
\newtheorem{proposition}[theorem]{Proposition}
\newtheorem{corollary}[theorem]{Corollary}
\newtheorem{conjecture}[theorem]{Conjecture}
\newtheorem{lemma}[theorem]{Lemma}

\newtheorem{remark}[theorem]{Remark}

\newcommand{\bt}{\begin{theorem}}
\newcommand{\et}{\end{theorem}}
\newcommand{\bl}{\begin{lemma}}
\newcommand{\el}{\end{lemma}}
\newcommand{\bp}{\begin{proposition}}
\newcommand{\ep}{\end{proposition}}
\newcommand{\bcor}{\begin{corollary}}
\newcommand{\ecor}{\end{corollary}}
\newcommand{\br}{\begin{remark}\rm}
\newcommand{\er}{\end{remark}}
\newcommand{\bcon}{\begin{conjecture}}
\newcommand{\econ}{\end{conjecture}}

\newcommand{\be}{\begin{equation}}
\newcommand{\ee}{\end{equation}}
\newcommand{\ba}{\begin{array}}
\newcommand{\ea}{\end{array}}
\newcommand{\bc}{\be\begin{array}{r@{\,}c@{\,}l}}
\newcommand{\ec}{\end{array}\ee}

\newcommand{\de}{\delta}
\newcommand{\De}{\Delta}
\newcommand{\eps}{\varepsilon}

\newcommand{\Gi}{{\cal G}}

\newcommand{\N}{{\mathbb N}}
\newcommand{\Z}{{\mathbb Z}}

\renewcommand{\P}{{\mathbb P}}
\newcommand{\E}{{\mathbb E}}

\newcommand{\sub}{\subset}

\newcommand{\ti}{\tilde}

\newcommand{\ffrac}[2]{{\textstyle\frac{{#1}}{{#2}}}}

\newcommand{\di}{\mathrm{d}}
\newcommand{\half}{{[0,\infty)}}

\begin{document}

\makeatletter\@addtoreset{equation}{section}
\makeatother\def\theequation{\thesection.\arabic{equation}} 

\renewcommand{\labelenumi}{{(\roman{enumi})}}

\title{\vspace{-3cm}Tightness of voter model interfaces}
\author{
Anja Sturm\vspace{6pt}\\
{\small Department of Mathematical Sciences}\\
{\small University of Delaware}\\
{\small 501 Ewing Hall}\\
{\small Newark, DE 19716-2553, USA}\\
{\small e-mail: sturm@math.udel.edu}\vspace{8pt}
\and Jan M. Swart\vspace{6pt}\\
{\small \' UTIA}\\
{\small Pod vod\'arenskou v\v e\v z\' i 4}\\
{\small 18208 Praha 8}\\
{\small Czech Republic}\\
{\small e-mail: swart@utia.cas.cz}\vspace{4pt}}
\date{{\scriptsize\file\quad\today
}}
\maketitle\vspace{-30pt}

\begin{abstract}\noi
Consider a long-range, one-dimensional voter model started with all
zeroes on the negative integers and all ones on the positive
integers. If the process obtained by identifying states that are
translations of each other is positively recurrent, then it is said
that the voter model exhibits interface tightness. In 1995, Cox and
Durrett proved that one-dimensional voter models exhibit interface
tightness if their infection rates have a finite third
moment. Recently, Belhaouari, Mountford, and Valle have improved this
by showing that a finite second moment suffices. The present paper
gives a new short proof of this fact. We also prove interface
tightness for a long range swapping voter model, which has a mixture
of long range voter model and exclusion process dynamics.

\end{abstract}

\noi
{\it Acknowledgements.} Work sponsored by GA\v CR grants 201/06/1323
and 201/07/0237, and by UDRF grant 06000596.

\section{Introduction and main results}\label{S:intro}

Let $X=(X_t)_{t\geq 0}$ be a long-range, one-dimensional `swapping' voter
model, i.e. $X$ is a Markov process with state space $\{0,1\}^\Z$
and formal generator $G:=G^{\rm v}+G^{\rm s}$, where
\bc\label{GX}
\dis G^{\rm v}f(x)&:=&\dis
\sum_{ij}q(i-j)1_{\{x(i)\neq x(j)\}}\{f(x^{\{i\}})-f(x)\},\\[5pt]
\dis G^{\rm s}f(x)&:=&\dis
\ffrac{1}{2}\sum_{ij}p(i-j)1_{\{x(i)\neq x(j)\}} \{f(x^{\{i,j\}})-f(x)\},
\ec
$q$ and $p$ are functions $\Z\to\half$, and $p$ is symmetric, i.e.,
$p(i)=p(-i)$ $(i\in\Z)$. Here, for any $x\in\{0,1\}^\Z$ and
$\De\sub\Z$,
\be
x^\De(i):=\left\{\ba{ll}
1-x(i)\quad&\mbox{if }i\in\De,\\
x(i)\quad&\mbox{if }i\not\in\De
\ea\right.
\ee
denotes the configuration obtained from $x$ by flipping all spins in
$\De$. Note that (\ref{GX}) says that if $X_t(i)\neq X_t(j)$ for some
$i\neq j$, then due to the action of $G^{\rm v}$, $X_t(i)$ adopts the
type of $X_t(j)$ with {\em infection rate} $q(i-j)$. This describes a
long-range, one-dimensional voter model. In addition, due to the
action of $G^{\rm s}$, at {\em swapping rate} $p(i-j)$ the sites $i$
and $j$ swap their types. In order for the process to be well-defined
(see \cite{Lig85}), we assume that $\sum_i(q(i)+p(i))<\infty$. We will
also require irreducibility conditions. We call a rate function
$r:\Z\to\half$ {\em irreducible} if each $i\in\Z$ can be written as
$i=i_1+\cdots+i_n$ with $n\geq 0$ and $r(i_k)>0$ for all $1\leq k\leq
n$.

Consider the space
\be
S_{\rm int}:=\big\{x\in\{0,1\}^\Z:
\lim_{i\to-\infty}x(i)=0,\ \lim_{i\to\infty}x(i)=1\big\}
\ee
of states describing the interface between two infinite population of
zeroes and ones. If the infection and swapping rates have a finite first
moment, i.e., $\sum_i|i|(q(i)+ p(i))<\infty$, then $X_0\in S_{\rm
int}$ implies $X_t\in S_{\rm int}$ for all $t\geq 0$, a.s.\ (see
\cite{BMV07} for this statement concerning $G^{\rm v}$).

Define an equivalence relation on $S_{\rm int}$ by setting $x\sim y$
if $x$ and $y$ are translations of each other and let $\ti S_{\rm
int}:=\{\ti x:x\in S_{\rm int}\}$, with $\ti x:=\{y\in S_{\rm
int}:y\sim x\}$, denote the set of equivalence classes. Then $\ti
X=(\ti X_t)_{t\geq 0}$ is a continuous-time Markov process with
countable state space $\ti S_{\rm int}$. Let $x_{\rm H}(i):=1_{\{i\geq
0\}}$ denote the `Heaviside state'. As long as $q(i)>0$ for some $i\neq 0$,
it is not hard to see that $\ti x_{\rm H}$ can be reached from any
state in $\ti S_{\rm int}$, hence $\ti X$ is an irreducible Markov
process on the set of all states that can be reached from $\ti x_{\rm
H}$. Following Cox and Durrett \cite{CD95}, we say that $X$ exhibits
{\em interface tightness} if $\ti X$ is positively recurrent on this
set.

Theorem~4 in \cite{CD95} states that long-range voter models (without
swapping) exhibit interface tightness provided that their infection
rates $q$ are symmetric, irreducible, and have a finite third
moment. Although not carried out there, their proofs also work if $q$
is asymmetric and $q_{\rm s}$ is irreducible, where $q_{\rm
s}(i):=\ffrac{1}{2}(q(i)+q(-i))$ denote the {\em symmetrized infection
rates}. Recently, Belhaouari, Mountford, and Valle \cite{BMV07}
improved this result by showing that a finite second moment
suffices. They showed that this condition is sharp in the sense that
interface tightness does not hold if $\sum_i|i|^cq(i)=\infty$ for some
$c<2$. We will give a new, short proof of their sufficiency result. In
fact, we will prove more:
\bt{\bf(Interface tightness for long-range swapping voter models)}
\label{T:long}\\
Assume that $\sum_i|i|^2(q(i)+p(i))<\infty$, $q(i)>0$ for some $i\neq 0$,
and that $q_{\rm s}+p$ is irreducible. Then $X$ exhibits interface tightness.
\et
For the symmetric nearest-neighbor case, Theorem~\ref{T:long} has been
proved in \cite[Theorem~1.1~(ii)~(b)]{BFMP01}. They also have results
for asymmetric exclusion dynamics (but symmetric voter dynamics). We
note that the symmetric nearest-neighbor swapping voter model arises
as the dual and interface model of certain systems of parity
preserving branching and annihilating random walks, see
\cite[Section~2.1]{SS07} and references there.

\section{Proof of Theorem~\ref{T:long}}

The main idea of the proof is the same as in \cite{CD95}, namely, to
look at the number of `inversions':
\be\label{fCD}
f_{\rm CD}(x):=|\{(i,j)\in\Z^2,\ i<j,\ x(i)>x(j)\}|\qquad(x\in S_{\rm int}).
\ee
In \cite[Section~4]{CD95} and \cite{BMV07}, this quantity is estimated
using duality and (subtle) results about one-dimensional random
walk. Our approach will be to insert $f_{\rm CD}$ into the generator
$G$ and prove that if interface tightness would not hold, then $f_{\rm
CD}$ would decrease unboundedly as time tends to infinity, which
yields a contradiction. In a way, our proof is similar to the methods
in \cite{BFMP01}, which are based on Lyapunov functions. The function
$f_{\rm CD}$, however, is not a Lyapunov function for our processs, so
our proofs need a probabilistic ingredient as well, which is provided
by Proposition~\ref{P:Igrow} below.

We start by calculating $Gf_{\rm CD}$.
\bl{\bf(Changes in number of inversions)}\label{L:change}
For each $x\in S_{\rm int}$, one has
\be\label{change}
Gf_{\rm CD}(x)=\sum_{n=1}^\infty\big(q_{\rm s}(n)+p(n)\big)n^2
-\sum_{n=1}^\infty q_{\rm s}(n)I_n(x),
\ee
where
\be\label{interface-n}
I_n(x):=|\{i\in\Z:x(i)\neq x(i+n)\}|\qquad(n\geq 1).
\ee
\el
We will need  the following lemma:
\bl{\bf(Nonnegative expectation)}\label{L:NE}
Assume that $\sum_i|i|^2(q(i)+p(i))<\infty$. Then
\be\label{NE}
\E[f_{\rm CD}(X_0)]+\int_0^t\E\big[Gf_{\rm CD}(X_s)\big]\,\di s\geq 0\qquad(t\geq 0).
\ee
\el
We will also need a result stating that the number of `boundaries'
between zeroes and ones, defined in the sense of (\ref{interface-n}),
grows over time if interface tightness does not hold.
\bp{\bf(Interface growth)}\label{P:Igrow}
Assume that $\sum_i|i|(q(i)+p(i))<\infty$, $q(i)>0$ for some $i\neq
0$, and $q_{\rm s}+p$ is irreducible. Assume that interface tightness
does not hold. Then
\be\label{Igrown}
\lim_{T\to\infty}\frac{1}{T}\int_0^T\!\!\di t\;\P[I_n(X_t)<N]=0
\qquad(N,n\geq 1).
\ee
\ep
With these statements we are now in the position to prove
Theorem~\ref{T:long}. Assume that interface tightness does not
hold. By our assumptions on $q$ and $p$ we can choose $i,N\geq 1$ such
that
\be\label{constant}
\sum_{n=1}^\infty(q_{\rm s}(n)+p(n))n^2<q_{\rm s}(i)N.
\ee
Then, by Lemmas \ref{L:NE} and \ref{L:change}, and Proposition~\ref{P:Igrow},
the process started in $X_0=x_{\rm H}$ satisfies
\begin{eqnarray}\label{minplus2}
\dis0&\leq&
\int_0^T\!\!\di t\;\E\big[Gf_{\rm CD}(X_t)\big]\\
\nonumber
&=&T\sum_{n=1}^\infty(q_{\rm s}(n)+p(n))n^2
-\sum_{n=1}^\infty q_{\rm s}(n)\int_0^T\!\!\di t\;\E\big[I_n(X_t)\big]\\
\nonumber
&\leq&T\sum_{n=1}^\infty(q_{\rm s}(n)+p(n))n^2
-q_{\rm s}(i)N\int_0^T\!\!\di t\;\P\big[I_i(X_t)\geq N\big]\\
\nonumber
&=&T\sum_{n=1}^\infty(q_{\rm s}(n)+p(n))n^2-(T-o(T))q_{\rm s}(i)N,
\end{eqnarray}
as $T\to\infty$. Due to (\ref{constant}) the right-hand side of
(\ref{minplus2}) tends to $-\infty$ as $T\to\infty$, which yields a
contradiction.\qed

\section{Proof of Lemmas~\ref{L:change} and \ref{L:NE}, and
Proposition~\ref{P:Igrow}}

{\bf Proof of Lemma~\ref{L:change}} We need to count the number of
pairs of sites $i,j$ with $i<j$ and $x(i)>x(j)$ that are created and
deleted due to the various possible jumps. We will first consider
$G^{\rm s} f_{\rm CD}$, i.e. the changes due to swapping. So consider
the case that $x(i)\neq x(i+n)$ for some $n \in \N$, while there are
$l$ ones on the left of $i$, $r$ zeroes on the right of $i+n$, and
$n_0$ zeroes and $n_1$ ones between $i$ and $i+n$. Then the changes in
$f_{\rm CD}(x)$ due to swapping can be summarized as follows:
\be\ba{cclr}\label{fCDflip2}
\underbrace{\ldots}_{l\times 1}\quad 0
\underbrace{\ldots}_{n_0\times 0,\ n_1\times 1}1\quad
\underbrace{\ldots}_{r\times 0}&\to&
\ldots 1\ldots 0\ldots &\phantom{-(}n_0+n_1+1\phantom{)}=\phantom{-}n\\[20pt]
\underbrace{\ldots}_{l\times 1}\quad1
\underbrace{\ldots}_{n_0\times 0,\ n_1\times 1}0\quad
\underbrace{\ldots}_{r\times 0}&\to&
\ldots 0\ldots 1\ldots &-(n_0+n_1+1)=-n\\
\ea
\ee
Thus, if we define
\be
I^{ab}_n(x):=|\{i:x(i)=a,\ x(i+n)=b\}|\qquad(n\geq 0,\ ab=01,10)
\ee
then we obtain
\be
\label{Gsf_CD1}
G^{\rm s}f_{\rm CD}(x)=\sum_{n=1}^\infty p(n)(n\cdot I^{01}_n(x)-n\cdot I^{10}_n(x))
\ee
Now, for any $0\leq m<n$, set
\be
I^{ab}_{n,m}(x):=|\{r\in\Z:x(nr+m)=a,\ x(n(r+1)+m)=b\}|.
\ee
Walking along the thinned-out lattice $n\Z+m$ from $-\infty$ to
$+\infty$, we see one more change from $0$ to $1$ than we see changes
{f}rom $1$ to $0$, i.e., $I_{n,m}^{01}(x)=I_{n,m}^{10}(x)+1$ for all
$x\in S_{\rm int}$. Since
$I^{ab}_n(x)=\sum_{m=0}^{n-1}I^{ab}_{n,m}(x)$, it follows that
\be
\label{In^01In^10}
I_n^{01}(x)=I_n^{10}(x)+n\qquad(x\in S_{\rm int}).
\ee
This implies that (\ref{Gsf_CD1}) simplifies to
\be
\label{Gsf_CD2}
G^{\rm s}f_{\rm CD}(x)=\sum_{n=1}^\infty p(n)n^2.
\ee
In order to consider the effect of $G^{\rm v}$ on $f_{\rm CD}$ we write
\bc\label{Gf}
\dis G^{\rm v}f_{\rm CD}(x)&=&\dis\sum_kq(-k)
\Big(\sum_{i:\,x(i)=1}\sum_{j:\,j>i}\big(1_{\{x(j)=1,\ x(j+k)=0\}}
-1_{\{x(j)=0,\ x(j+k)=1\}}\big)\\[5pt]
&&\dis+\sum_{i:\,x(i)=0}\sum_{j:\,j<i}\big(1_{\{x(j)=0,\ x(j+k)=1\}}
-1_{\{x(j)=1,\ x(j+k)=0\}}\big)\Big).
\ec
We observe that for any $x\in S_{\rm int}$ and $k>0$,
\be\ba{l}\label{cancel1}
\dis\sum_{j:\,j>i}\big(1_{\{x(j)=1,\ x(j+k)=0\}}
-1_{\{x(j)=0,\ x(j+k)=1\}}\big)\\[5pt]
\quad\dis=\sum_{j=i+1}^{i+k}\sum_{n\geq 0}
\big(1_{\{x(j+nk)=1,\ x(j+(n+1)k)=0\}}
-1_{\{x(j+nk)=0,\ x(j+(n+1)k)=1\}}\big)\\[5pt]
\quad\dis=-\sum_{j=i+1}^{i+k}1_{\{x(j)=0\}}.
\ec
To see why the last equality in (\ref{cancel1}) holds, observe that
since $k>0$, the sequence $x(j),x(j+k),x(j+2k)\ldots$ is eventually
one. Hence, if $x(j)=1$, then the number of changes from $0$ to $1$
equals the number of changes from $1$ to $0$ and all terms cancel,
while if $x(j)=0$ there is one extra change from $0$ to $1$, leading
to a contribution of minus one.

Likewise, for any $x\in S_{\rm int}$ and $k>0$,
\be\ba{l}\label{cancel2}
\dis\sum_{j:\,j<i}\big(1_{\{x(j)=0,\ x(j+k)=1\}}
-1_{\{x(j)=1,\ x(j+k)=0\}}\big)\\[5pt]
\quad\dis=\sum_{j=i-k}^{i-1}\sum_{n\geq 0}
\big(1_{\{x(j-nk)=0,\ x(j-(n-1)k)=1\}}
-1_{\{x(j-nk)=1,\ x(j-(n-1)k)=0\}}\big)\\[5pt]
\quad\dis=\sum_{j=i}^{i+k-1}\sum_{n\geq 0}
\big(1_{\{x(j-(n+1)k)=0,\ x(j-nk)=1\}}
-1_{\{x(j-(n+1)k)=1,\ x(j-nk)=0\}}\big)\\[5pt]
\quad\dis=\sum_{j=i}^{i+k-1}1_{\{x(j)=1\}}.
\ec
It follows that for any $k>0$,
\be\ba{l}\label{calc}
\dis\sum_{i:\,x(i)=1}\sum_{j:\,j>i}\big(1_{\{x(j)=1,\ x(j+k)=0\}}
-1_{\{x(j)=0,\ x(j+k)=1\}}\big)\\[5pt]
\dis+\sum_{i:\,x(i)=0}\sum_{j:\,j<i}\big(1_{\{x(j)=0,\ x(j+k)=1\}}
-1_{\{x(j)=1,\ x(j+k)=0\}}\big)\\[5pt]
\dis\quad=-\sum_i\sum_{j=i+1}^{i+k}1_{\{x(i)=1,\ x(j)=0\}}
+\sum_i\sum_{j=i}^{i+k-1}1_{\{x(i)=0,\ x(j)=1\}}\\[5pt]
\dis\quad=\sum_{n=1}^{k-1}I^{01}_n(x)-\sum_{n=1}^kI^{10}_n(x),
\ec
Using (\ref{In^01In^10}), the expression in (\ref{calc}) can be rewritten as
\be
-I^{10}_k(x)+\sum_{n=1}^{k-1}(I^{01}_n(x)-I^{10}_n(x))=
-\ffrac{1}{2}(I_k(x)-k)+\sum_{n=1}^{k-1}n=\ffrac{1}{2}(k^2-I_k(x)),
\ee
which holds for $k>0$. Using symmetry with respect to the map
$x\mapsto x'$ where $x'(i):=1-x(-i)$, it is not hard to see that we
get the same formula for the expression in (\ref{calc}) if
$k<0$. Inserting this into (\ref{Gf}), we arrive at 
\be
\label{Gvf_CD}
G^{\rm v} f_{\rm CD}(x)= \sum_{n=1}^\infty q_{\rm s}(n) (n^2-I_n(x)).
\ee
Taking (\ref{Gsf_CD2}) and  (\ref{Gvf_CD}) together now implies
(\ref{change}).\qed

\noi
{\bf Proof of Lemma~\ref{L:NE}} If the function $f_{\rm CD}$ were
bounded, then standard theory would tell us that the process
\be\label{MP}
M_t:=f_{\rm CD}(X_t)-\int_0^t Gf_{\rm CD}(X_s)\di s\qquad(t\geq 0)
\ee
is a martingale with respect to the filtration generated by $X$. In
particular, since $\E[M_t]=\E[M_0]$ and $f_{\rm CD}\geq 0$, this would
imply (\ref{NE}). In the present case, since $f_{\rm CD}$ is
unbounded, we have to work a bit. Let
\be\label{wdef}
w(x):=\max\{i:x(i)\neq x(i+1)\}-\min\{i:x(i)\neq x(i+1)\}
\ee
denote the `width' of an interface state $x\in S_{\rm int}$. We
can couple $X_t$ to a continuous-time random walk $(R_t)_{t\geq 0}$,
started in $R_0=w(X_0)$, which jumps from $r$ to $r+n$ with rate
\be\label{adef}
a(n):=\sum_{k=n}^\infty2(q_{\rm s}(k)+p(k)),
\ee
in such a way that
$R_t\geq w(X_t)$ for all $t\geq 0$ a.s. We check that
\be\ba{l}\label{a1}
\dis\sum_{n=1}^\infty a(n)n
=\sum_{n=1}^\infty\sum_{k=n}^\infty2(q_{\rm s}(k)+p(k))n
=\sum_{k=1}^\infty\sum_{n=1}^k2(q_{\rm s}(k)+p(k))n\\[5pt]
\quad\dis=\sum_{k=1}^\infty2(q_{\rm s}(k)+p(k))\ffrac{1}{2}k(k+1)<\infty.
\ec
Let $X^K$ be the process with generator given by
\bc\label{GK}
\dis G^Kf(x)&:=&\dis
\!\!\sum_{ij:\;|i-j|\leq K}\!\!
q(i-j)1_{\{x(i)\neq x(j)\}}\{f(x^{\{i\}})-f(x)\},\\[5pt]
&&\dis+\ffrac{1}{2}\!\!\sum_{ij:\;|i-j|\leq K}\!\!
p(i-j)1_{\{x(i)\neq x(j)\}} \{f(x^{\{i,j\}})-f(x)\},
\ec
started in $X^K_0=X_0$. We set
\be\label{taudef}
\tau_{K,N}:=\inf\{t\geq 0:w(X^K_t)>N\}
\ee
and define $\tau_N$ similarly, with $X^K$ replaced by $X$. It follows
{f}rom standard theory that
\be\label{MPKN}
M^{K,N}_t:=f_{\rm CD}(X^K_{t\wedge\tau_{K,N}})
-\int_0^{t\wedge\tau_{K,N}}\!\!G^Kf_{\rm CD}(X^K_s)\,\di s\qquad(t\geq 0)
\ee
is a martingale, where by Lemma~\ref{L:change},
\be\label{Kchange}
G^Kf_{\rm CD}(x)=\sum_{n=1}^K\big(q_{\rm s}(n)+p(n)\big)n^2
-\sum_{n=1}^Kq_{\rm s}(n)I_n(x).
\ee
Using an obvious coupling, we may arrange that $X^K$ converges a.s.\
to $X$ in such a way that there exists a random $K_0$ such that for
all $K\geq K_0$, one has $\tau_{K,N}=\tau_N$ and $X^K_t=X_t$ for all
$t\leq\tau_N$. Since moreover $\tau_N\to\infty$ as $N\to\infty$ (which
follows from the estimate $w(X_t)\leq R_t$), taking into account
(\ref{change}) and (\ref{Kchange}), it follows that
\be\label{as}
\lim_{N\to\infty}\lim_{K\to\infty}\int_0^{t\wedge\tau_{K,N}}\!\!
G^Kf_{\rm CD}(X^K_s)\,\di s=\int_0^tGf_{\rm CD}(X_s)\,\di s\quad{\rm a.s.}
\ee
Since $M^{K,N}$ is a martingale,
\be\label{preNE}
0\leq\E\big[f_{\rm CD}(X^K_{t\wedge\tau_{K,N}})\big]
=\E\big[f_{\rm CD}(X_0)\big]
+\E\Big[\int_0^{t\wedge\tau_{K,N}}\!\!G^Kf_{\rm CD}(X^K_s)\,\di s\Big].
\ee
Letting $K\to\infty$ and then $N\to\infty$ in (\ref{preNE}), using
(\ref{as}), we arrive at (\ref{NE}), provided we show (in view of
(\ref{Kchange})) that the random variables
\be
\sum_{n=1}^Kq_{\rm s}(n)\int_0^{t\wedge\tau_{K,N}}\!\!I_n(X^K_s)\,\di s
\qquad(K,N\geq 1)
\ee
are uniformly integrable. We can couple the $X^K$ to a random walk
$(R_t)_{t\geq 0}$ with jump rates $a(n)$ as in (\ref{adef}), in such a
way that $R_t\geq w(X^K_t)$ for all $t\geq 0$ and $K\geq 1$
a.s. Hence, since $I_n(x)\leq n+w(x)$, we may estimate
\be\label{ui}
\sum_{n=1}^Kq_{\rm s}(n)
\int_0^{t\wedge\tau_{K,N}}\!\!I_n(X^K_s)\,\di s
\leq\sum_{n=1}^\infty q_{\rm s}(n)\int_0^t(n+R_s)\,\di s.
\ee
By (\ref{a1}), the right-hand side of (\ref{ui}) has finite
expectation, proving the required uniform integrability.\qed

\br
{\bf(Martingale problem)}
If $\sum_i|i|^3(q(i)+p(i))<\infty$, then the jump rates $a(n)$ in
(\ref{adef}) have a finite second moment. Using this and the estimate
$f_{\rm CD}(x)\leq w(x)^2$, one can prove that in this case the
process in (\ref{MP}) is a martingale.
\er
{\bf Proof of Proposition~\ref{P:Igrow}} Consider the `boundary process'
\be\label{bp}
Y_t(i):=1_{\{X_t(i)\neq X_t(i+1)\}}\qquad(t\geq 0,\ i\in\Z),
\ee
which is a Markov process in 
\be\label{Sb}
S_{\rm bound}:=\big\{y\in\{0,1\}^\Z:\sum_iy(i)\mbox{ is finite and odd}\big\}.
\ee
Then $(Y_t)_{t\geq 0}$ is a Markov process in $\{0,1\}^\Z$
with formal generator
\bc
\label{Ygen}
\dis G_Yf(y)&:=&\dis\sum_{i<j}1_{\txt\{\sum_{k=i+1}^jy(k)\mbox{ is odd}\}}
q(i-j)\big\{f(y^{\{i,i+1\}})-f(y)\big\}\\[5pt]
&&\dis+\sum_{i>j}1_{\txt\{\sum_{k=j+1}^iy(k)\mbox{ is odd}\}}
q(i-j)\big\{f(y^{\{i,i+1\}})-f(y)\big\}\\[5pt]
&&\dis+\sum_{i<j-1}1_{\txt\{\sum_{k=i+1}^jy(k)\mbox{ is odd}\}}
p(i-j)\big\{f(y^{\{i,i+1,j,j+1\}})-f(y)\big\}\\[5pt]
&&\dis+\sum_i1_{\txt\{y(i+1)=1\}}
p(1)\big\{f(y^{\{i,i+2\}})-f(y)\big\}.
\ec
Using the fact that $\sum_i|i|(q(i)+p(i))<\infty$, one can check that
$Y$ is a parity preserving cancellative spin system in the sense of
\cite{SS07}. We will show that the process started in $y$, denoted by
$Y^y$, satisfies
\be\label{extpos}
\inf\big\{\P\big[|Y^y_t|=n\big]:|y|=n+2,\ y(i)=1=y(j)
\mbox{ for some }i\neq j,\ |i-j|\leq L\big\}>0
\ee
$(L\geq 1,\ n\geq 0,\ t>0)$. As a result of (\ref{extpos}), we can
apply \cite[Proposition~13]{SS07}
to conclude that
\be\label{Igrow3}
\lim_{T\to\infty}\frac{1}{T}\int_0^T\!\!\di t\;\P[I_1(X_t)<N]=0\qquad(N\geq 1).
\ee
See also Proposition~2.6 in \cite{Han99} for similar arguments
concerning a dual process to the threshold voter model.

To give a rough idea of the proof and of the importance of condition
(\ref{extpos}), think of the sites with $Y_t(i)=1$ as being occupied
by a particle. Then $Y$ is a parity preserving particle system, i.e.,
if $Y$ is started in an odd (even) initial state, then the number of
particles always stays odd (even).  Let $\ti Y$ denote the process
obtained from $Y$ by identifying states that are a translation of each
other. If $\ti X$ is not positively recurrent, then the same is true
for $\ti Y$. The main idea of the proof of (\ref{Igrow3}) is to use
induction on $n$ to show that there cannot be less than $n$ particles
for a positive fraction of time. This is obviously true for $n=1$;
imagine that it holds for a certain $n$. If we see $n$ particles for a
positive fraction of time, then most of the time these particles must
be situated far from each other, for else with positive probability at
least two would annihilate each other due to (\ref{extpos}), violating
the induction hypothesis. However, since $\ti Y$ is not positively
recurrent, $n$ single particles far from each other will soon each
produce three particles, hence we cannot see $n$ particles for a
positive fraction of time.

To boost up (\ref{Igrow3}) to the statement in (\ref{Igrown}), it
suffices to show that the process started in $x$, denoted by $X^x$,
satisfies
\be\label{boost}
\lim_{N\to\infty}\inf_{|I_1(x)|\geq N}\P[I_n(X^x_t)<M]=0
\qquad(M,n\geq 1,\ t>0).
\ee
For if (\ref{boost}) holds, then for each $t,\eps>0$ we can choose $N$
large enough such that the limit in (\ref{boost}) is smaller than $\eps$,
and therefore, by (\ref{Igrow3}) and a restart argument
\be
\limsup_{T\to\infty}\frac{1}{T}\int_t^{t+T}\!\!\di s\;\P[I_n(X_s)<M]
\leq\eps\qquad(M\geq 1).
\ee
Since $\eps>0$ is arbitrary, this implies (\ref{Igrown}).

We still need to prove (\ref{extpos}) and (\ref{boost}). We start with
the former. Choose $k\neq 0$ such that $q(k)>0$. By symmetry, we may
without loss of generality assume that $k>0$. It suffices to prove
(\ref{extpos}) for $L\geq k$. So fix $L\geq k$, $n\geq 0$, and
$t>0$. Let $y=\de_{i_1}+\cdots+\de_{i_n}$ where $\de_i(j):=1$ if $i=j$
and $=0$ otherwise, and $i_1<\cdots<i_n$. Assume that
\be
M:=\big\{m'\in\{1,\ldots,n-1\}:i_{m'+1}-i_{m'}\leq L\big\}
\ee
is not empty and let $m:=\inf(M)$. Then
$x(i_m+1-L)=\cdots=x(i_m)\neq x(i_m+1)$, and hence, since $L\geq k$,
there is a positive probability that during the time interval $[0,t]$,
the sites $i_m+1,\ldots,i_{m+1}$  get infected successively by the sites
$i_m+1-k,\ldots,i_{m+1}-k$, leading to a decrease in $I_1(X_t)$ of
2. Using the fact that $\sum_i|i|(q(i)+p(i))<\infty$, it is not hard
to see that moreover, with positive probability, no other infections
take place, and that this probability is uniformly bounded from below in
all $y$ satisfying our assumptions.

To prove (\ref{boost}), we view the dynamics of our process $X$ as
follows. For each ordered pair $(i,j)$ with $i\neq j$, at times
selected according to an independent Poisson point process with
intensity $q(i-j)$, the type of site $j$ infects the site
$i$. Likewise, for each unordered pair $\{i,j\}$ with $i\neq j$, at
times selected according to an independent Poisson point process with
intensity $p(i-j)$, the sites $i$ and $j$ swap their types.

We claim that if we view the evolution of types in this way, then we
can find a $j\in\Z$ such that with positive probability $X_t(j)$
inherits its type from $x(0)$ and $X_t(j+n)$ inherits its type from
$x(1)$. To see this, we will make a number of infections and swaps to
transport the type of site $0$ to $j$ and the type of site $1$ to
$j+n$. Let us say that in the $k$-th step of our construction, we have
transported the type of site $0$ to the site $l_k$ and the type of
site $1$ to $r_k$. Then, in the $(k+1)$-th step of our construction,
by making an infection with rate $q(i)$, we may transport the type of
site $l_k$ to $l_k+i$ and set $(l_{k+1},r_{k+1}):=(l_k+i,r_k)$, or we
may transport the type of site $r_k$ to $r_k+i$ and set
$(l_{k+1},r_{k+1}):=(l_k,r_k+i)$, and similarly for swaps. Thus, in
each step, we may increase $r_k-l_k$ by $d$ for each
$d\in\Gi:=\{i\in\Z:q_{\rm s}(i)+p(i)>0\}$. We must only make sure that
when we move the type of one site to a new position, we do not
influence the type of the other site. To avoid this sort of influence,
we will make sure that at each point in our construction $l_k<r_k$. We
claim that this is possible. Set $n_k:=r_k-l_k$. We need to show that
there exist {\em (strictly) positive} integers $n_0,\ldots,n_m$ such that
$n_0=1$, $n_m=n$, and $(n_k-n_{k-1})\in\Gi$ for all $1\leq k\leq
m$. By our assumption that $q_{\rm s}+p$ is irreducible we can write
$n-1=i_1+\cdots+i_m$ with $i_1,\ldots,i_m\in\Gi$. (Note that this is
the only place in our proofs where we use irreducibility.) Without
loss of generality we may assume that $i_1\geq\cdots\geq i_m$. Then
setting $n_k:=1+i_1+\cdots+i_k$ proves our claim.

It follows that there exist $j\in\Z$ and $L\geq j$ such that whenever
$x(i)\neq x(i+1)$, there is a positive probability that $X_t(i+j)$
inherits its type from $x(i)$ and $X_t(i+j+n)$ inherits its type from
$x(1)$ through a sequence of infections and swaps that are entirely
contained in $\{i-L,\ldots,i+L\}$. If $I_1(x)$ is large, we can find
many sites $i$, situated at least a distance $3L$ from each other,
such that $x(i)\neq x(i+1)$. By what we have just proved each pair
has an independent probability to produce at time $t$ two sites $i+j$
and $i+j+n$ such that $X_t(i+j)\neq X_t(i+j+n)$, hence $I_n(X_t)$ is
with large probability large.\med

\noi
{\bf\large Acknowledgement} We thank Rongfeng Sun who found an error
in our original statement of Lemma~\ref{L:NE} and told us how to fix
it.

\newcommand{\noopsort}[1]{}

\end{document}